\documentclass[12pt]{article}

\usepackage{amssymb}

\newtheorem {theorem} {Theorem}

\newtheorem {lemma}  [theorem]{Lemma}
\newtheorem {remark} [theorem]{\sc Remark}

\newcommand{\bbox}{\ \hfill\rule[-1mm]{2mm}{3.2mm}}

\title{On the cyclicity of weight--homogeneous centers \thanks{The second and third authors are partially supported by a MCYT/FEDER
grant number MTM2008-00694. The second author is also partially
supported by a CIRIT grant number 2005SGR 00550.}}

\author{{\sc Lubomir Gavrilov$^{\ (1)}$, Jaume Gin\'e$^{\ (2)}$} \\ {\sc and Maite Grau$^{\ (2)}$}}
\date{}

\begin{document}
\maketitle \vspace{-1.0cm}

\begin{center}
\noindent {\small{$^{\ (1)}$ Institut de Math\'{e}matiques, UMR 5219\\
 Universit\'{e} de Toulouse\\  31062 Toulouse, Cedex 9,  France   \\ {\rm E--mail:} {\tt
lubomir.gavrilov@math.ups-tlse.fr} \\
$^{\ (2)}$ Departament de Matem\`atica. Universitat de Lleida. \\
Avda. Jaume II, 69. 25001 Lleida, Spain \\ {\rm E--mails:} {\tt
gine@matematica.udl.cat},\ {\tt mtgrau@matematica.udl.cat}}}
\end{center}
\begin{abstract}
Let $W$ be a weight--homogeneous planar polynomial differential
system with a center. We find an upper bound of the number of
limit cycles which bifurcate from the period annulus of $W$  under
a generic polynomial perturbation. We apply this result to a
particular family  of planar polynomial systems  having a
nilpotent center without meromorphic first integral.
\end{abstract}

{\small{\noindent 2000 {\it AMS Subject Classification:} 34C07, 34C05, 34C14.  }}

\section{Statement of the results}

The real planar differential system
\begin{equation}
\label{1}
\dot{x}\, =\,
f(x,y), \quad \dot{y} \, = \, g(x,y) \label{eqg} \end{equation}
is said to be $(\alpha,\beta,\omega)$--{\em weight--homogeneous}  if there exist
weights $\alpha,\beta >0$ and $\omega\in\mathbb{R}$ such that for all $\rho
>0$, $x,y \in \mathbb{R}$ holds
$f(\rho^\alpha\, x, \rho^\beta\, y) \, = \, \rho^{\omega+\alpha} \, f(x,y)$ and
$g(\rho^\alpha\, x, \rho^\beta\, y) \, = \, \rho^{\omega+\beta} \, g(x,y)$. Equivalently, the system (\ref{1}) is weight--homogeneous
if the foliation
$$
\mathcal{F}_0 : f(x,y) \, dy - g(x,y)\, dx = 0
$$
is invariant under the dilatation
\begin{equation}
\label{dilatation}
\Phi : \mathbb{R}^2 \rightarrow \mathbb{R}^2 : (x,y) \mapsto \Phi(x,y)=(\rho^\alpha\, x, \rho^\beta\, y) .
\end{equation}

\par \emph{In what follows we shall suppose that the weight--homogeneous system {\rm
(\ref{eqg})} has a center at the origin}. Then the center is
global and the open period annulus $\mathcal{O} = \mathbb{R}^2
\setminus \{(0,0)\}$ is a union of periodic orbits. The center
needs to be global due to its invariance under the aforementioned
dilatation. The purpose of the present paper is to study small
perturbations of such  global centers.

Consider a one--parameter analytic perturbation
\begin{equation}
\dot{x} \, =\, f(x,y) + \, \varepsilon \, Q(x,y,\varepsilon), \quad \dot{y} \, =\, g(x,y) - \, \varepsilon \,
P(x,y,\varepsilon) \label{eqpa}
\end{equation}
of (\ref{eqg}) and suppose that $P(x,y,\varepsilon),
Q(x,y,\varepsilon)$ are polynomials in $x,y$ of degree $d$, depending analytically on the parameter $\varepsilon$.
For
every compact set $K$ contained in the real positive
half-axis $ \mathbb{R}^+_*=\{ (x,0): x>0\}$, and for $|\varepsilon|$
small enough, there exists an open interval $\Delta\supset K$ on which the first return map

\begin{equation}\label{return}
\Pi_\varepsilon: \Delta \rightarrow \mathbb{R}^+ : x \mapsto
\Pi_\varepsilon(x)
\end{equation}
is well defined and analytic
$$
\Pi_\varepsilon(x)= x + \varepsilon^k M_k(x)+ o(\varepsilon^k) .
$$
Its fixed points correspond to limit cycles of (\ref{eqpa}). We
note that the function $M_k(x)$ is defined on the whole half-axis $
\mathbb{R}^+_*=\{ (x,0): x>0\}$ and its zeros counted with
multiplicity provide an upper bound for the number of fixed points
of $\Pi_\varepsilon$ on $K$. It is known that the so called
higher order Poincar\'{e}-Pontryagin-Melnikov function $M_k(x)$
allows an integral representation in terms of iterated path
integrals \cite{Francoise,gav, Iliev} along the periodic orbits
$\gamma(x)$ of the system (\ref{eqg}). If the perturbation is
generic, then $k=1$ ($M_1 \not\equiv 0$). The first result of the
paper is the following
\begin{theorem}
\label{theo} If the first Poincar\'{e}-Pontryagin-Melnikov function $M_1 $ is not identically zero, then
 the perturbed system $(\ref{eqpa})$ has at most $(d+1)(d+4)/2-1$ limit cycles which tend to
periodic orbits as $\varepsilon$ tends to zero.
\end{theorem}
In other words, the cyclicity of the open period annulus $\mathcal{O}$ with respect
to generic (such that $M_1\not\equiv 0$) perturbations of degree
$d$ is at most $(d+1)(d+4)/2-1$. This  bound is certainly not
exact, as one can easily check in the case of a linear center. The above theorem does not make any claim
about the number of limit cycles which tend to the origin or to the "infinity".

 In
the case when the system (\ref{1}) has a polynomial first integral and is
moreover Hamiltonian, the computation is straightforward (as we
have a well known integral formula for $M_1$), see \cite{Li,
Llibre, ZhaoZhang}. However, in general, a weight--homogeneous
system with a global center is neither Hamiltonian, nor it has an
analytic or even meromorphic first integral. Under the
restrictions that the polynomials $f,g$ have no common divisor in
$\mathbb{R}[x,y]$ and that the origin has no characteristic
directions, a result close to our Theorem \ref{theo} was recently
announced in \cite[Theorems A, B, C]{LLYZ}. We note that our
result is different, more general and with a much shorter proof.

The next question addressed in the paper is: \emph{are there  non-Hamiltonian weight-homogeneous polynomial
systems with a center ?}

The answer turns out to be positive. Theorem \ref{prop3} provides
a large class of non-trivial weight-homogeneous systems with a
global center to which Theorem \ref{theo} applies. An explicit
example of such a system is given in (\ref{eq2+}) below, and the
exact upper bound for the number of the limit cycles under a
generic perturbation is given in
 Theorem \ref{th1}.

\par To formulate Theorem \ref{prop3}, let $h_1,h_2\in
\mathbb{R}[x,y]$ be two
$(\alpha,\beta,\delta)$--weight--homogeneous polynomials of the
same weighted degree, that is, $h_i(\rho^\alpha\, x, \rho^\beta\,
y) \, = \, \rho^{\delta} \, h_i(x,y)$, $i=1,2$, for any $x,y \in
\mathbb{R}$ and any $\rho \in \mathbb{R}_{+}$. We recall that
$\alpha, \beta, \delta$ are real numbers with $\alpha,\beta>0$. We
also assume that $h_1$ and $h_2$ are such that
\begin{itemize}
    \item $h_2(x,y)\, \geq\, 0$, $\forall (x,y)\in \mathbb{R}^2$
    \item $h_1(x,y)\, =\, h_2(x,y)\, =\, 0$ if and only if $(x,y)=(0,0)$ .
\end{itemize}
Let $\sigma \in \mathbb{C}$ and $\Re {\sigma} \neq 0$ and put
$$
H= (h_1+i h_2)^\sigma(h_1-i h_2)^{\bar{\sigma}}, \qquad V= \frac{1}{2}\, (h_1 + i h_2)^{1-\sigma}(h_1-i
h_2)^{1-\bar{\sigma}} .
$$
\begin{theorem}
The system
\begin{equation}\label{eq2}
\dot{x}\, = \, H_y\, V, \quad \dot{y}\, =\, -H_x\, V
\end{equation}
is a real polynomial
$(\alpha,\beta,2\delta-\alpha-\beta)$--weight--homogeneous planar
differential system which has a global center.
 \label{prop3}
\end{theorem}
The proof of this Theorem is given in Section \ref{sect2}. Clearly some hypothesis can be relaxed. For instance,
$h_1,h_2$ need not be polynomials.
\newline

{\bf Example} \emph{Put
$$
h_1(x,y)=x^{2n}+y, \quad  h_2(x,y)= |\sqrt{1+4c}|\,\, x^{2n}/2, \quad \sigma = 1- 1/\sqrt{1+4c},
$$
where $n$ is a natural number $n \geq 1$  and $c$ is a real number
with $c<-1/4$. System {\rm (\ref{eq2})} takes the form
\begin{equation}
\dot{x} \, = \, y \, + \, x^{2n}, \qquad \dot{y} \, = \, 2 \, n \, c \, \, x^{4n-1}, \label{eq1}
\end{equation}
and it has a global center with a first integral \begin{equation}
\label{Hn} H= (y+ \mu_{+} x^{2n})^{1- 1/\sqrt{1+4c}}(y+ \mu_{-}
x^{2n})^{1+ 1/\sqrt{1+4c}} , \end{equation} where $\mu_{\pm}=(1
\pm \sqrt{1+4c})/2$.} \newline

We remark that system (\ref{eq1}) is a
$(1,2n,2n-1)$--weight--homogeneous differential system. \newline

To apply Theorem \ref{theo} to (\ref{eq1}), we consider the
following perturbed system
\begin{equation}
\dot{x} \, = \, y \, + \, x^{2n} \, + \, \varepsilon \, Q(x,y,\varepsilon), \qquad \dot{y} \, = \, 2\, n\, c \,
\, x^{4n-1} \, - \, \varepsilon \, P(x,y,\varepsilon), \label{eq2+}
\end{equation}
where $|\varepsilon|>0$ is a small parameter, $P(x,y,\varepsilon)$ and $Q(x,y,\varepsilon)$ are polynomials in
$(x,y)$ and depend analytically on $\varepsilon$ and $P(x,y,0)$ and $Q(x,y,0)$ are polynomials of degree at most
$4n-1$. Let $\{\gamma_h\}_h$ be the continuous family of periodic orbits surrounding the center. The first
Poincar\'{e}-Pontryagin-Melnikov function $M_1$ can be written as follows
\begin{equation} \label{m1}
M_1(h) \, = \, \oint_{\gamma_h} \frac{ P(x,y,0) \, dx \, + \, Q(x,y,0) \, dy}{V(x,y)} \, . \label{eqi}
\end{equation}
(this formula holds true for the system (\ref{eqpa}) too, where $V$ is an appropriate integrating factor). We
deduce the following
\begin{theorem}
Suppose that  $M_1(h)$ is not identically zero. Then
\begin{enumerate}
\item[{\rm (a)}]
 the perturbed system $(\ref{eq2+})$ has at most $n(3n+1)-1$ limit cycles which tend to
period orbits as $\varepsilon$ tends to zero.
\item[{\rm (b)}] For suitable polynomials $P(x,y,0) ,Q(x,y,0) $, and for all sufficiently small $|\varepsilon|$,
 the perturbed system $(\ref{eq2+})$
has at least $ n(2n+1)-1$ limit cycles which tend to period orbits.
\end{enumerate} \label{th1}
\end{theorem}
The bounds are written in terms of the number $n$ associated to system (\ref{eq1}) which is not its degree. The
degree $d=4n-1$, and hence the upper bound for the number of limit cycles  is $(d+1)(3d+7)/16-1$, which is
better than the one given in Theorem \ref{theo}, due to the symmetries of (\ref{eq1}). The lower bound can also
be written in terms of the degree and it is $(d+1)(d+3)/8-1$. It is obtained by a direct study of the function
$M_1$. The symmetries of system (\ref{eq1}) imply that the lower bound for the number of limit cycles
($(d+1)(d+3)/8-1$) is strictly lower than the corresponding lower bound ($d(d+1)/2$) for a degree $d$
perturbation of a generic Hamiltonian system of degree $d$ \cite{Ilya}.

\section{Proofs. \label{sect2}}

{\em  Proof of Theorem} \ref{theo}. We shall estimate the number
of isolated zeros (counted with multiplicity) of the first
Poincar\'{e}-Pontryagin-Melnikov function $I(h)=M_1(h)$, see
(\ref{m1}). Consider first the special case $Q(x,y,0)
= 0$ and $P(x,y,0) = x^i y^j$ with $i,j \geq 0$ and $i+j\leq d$. Denote by $\cal{F}_\varepsilon$ the real foliation
defined by
\begin{equation}
\label{1eps}
{\cal F}_\varepsilon :
-g(x,y) \,dx +
f(x,y) \,dy + \varepsilon x^i y^j \,dx = 0
\end{equation}
and let $ \Pi_{\varepsilon}$ be the corresponding first return map.
The dilatation $\Phi$ defined in (\ref{dilatation})  transforms ${\cal F}_\varepsilon$ to the foliation
$\Phi^*{\cal F}_\varepsilon$
\begin{equation}
\label{11eps} \Phi^* {\cal F}_\varepsilon : -g(x,y) \,dx + f(x,y)
\,dy +  \varepsilon \rho^{\ell} x^i y^j \,dx = 0, \quad \ell=
\omega+\beta-i\alpha-j\beta.
\end{equation}
The first return map $\Pi_{\varepsilon \,
\rho^\ell}$ of the foliation  ${\cal F}_{\varepsilon \,
\rho^\ell}= \Phi^*{\cal F}_\varepsilon$ is therefore conjugated to $\Pi_{\varepsilon}$
\begin{equation}
\label{conjugated}
\Pi_{\varepsilon} = \phi^{-1}\circ\Pi_{\varepsilon \,
\rho^\ell} \circ \phi
\end{equation}
where $\phi(x)=\Phi(x,0)$ is the restriction of $\Phi$ to the
cross-section $ \mathbb{R}^+_*=\{ (x,0): x>0\}$. \par Given a
foliation $\cal{F}$ with an associated Poincar\'e map $\Pi$ and a
diffeomorphism $\Phi$, then the Poincar\'e map associated to
$\Phi^*{\cal F}$ is conjugated to $\Pi$. This is due to the fact
that the Poincar\'e map of a foliation is defined by its flow, see
for instance Theorem 1 in page 207 of \cite{Perko}, and the flows
of the foliations $\cal{F}$ and $\Phi^*{\cal F}$ are conjugated,
see for instance Lemma 11 in page 217 of \cite{Spivak}. \par We
have
\begin{eqnarray}
\Pi_{\varepsilon}(x)&= &x + \varepsilon \, I(x) + o(\varepsilon) \label{geq1}\\
\Pi_{\varepsilon\,\rho^\ell}(x)&=& x + \varepsilon \, \rho^\ell\, I(x) + o(\varepsilon)\\
\phi^{-1}\circ\Pi_{\varepsilon \,
\rho^\ell} \circ \phi(x) &= &\rho^{-\alpha }\,( \rho^{\alpha }\,x +  \varepsilon \, \rho^\ell \, I( \rho^{\alpha }\,x) + o(\varepsilon))\\
&=& x+  \varepsilon \, \rho^{\ell-\alpha }\, I( \rho^{\alpha
}\,x)+ o(\varepsilon) . \label{geq2}
\end{eqnarray}
Therefore, equating the first order terms in $\varepsilon$ in
(\ref{geq1}) and (\ref{geq2})
 we get
$$
I(x) \ = \ \rho^{\ell -\alpha}\, I(\rho^\alpha x),
$$
for any positive real numbers $\rho,x$. This implies, choosing
$\rho\, =\, x^{-1/\alpha}$, that
\begin{equation}
I(x) \, = \, I(1)
\, x^{(\alpha-\ell)/\alpha}= x^{1-(\omega+\beta)/\alpha} \, I(1) \, x^{(\alpha\,i+\beta
j)/\alpha}
\end{equation}

 In a similar way, if we suppose that $Q(x,y,0) = x^i y^j$
and $P(x,y,0) = 0$ with $i,j \geq 0$ and $i+j\leq d$,  we
deduce \begin{equation}
I(x) \, = \,
x^{1-(\omega+\beta)/\alpha} \, I(1) \, x^{(\alpha\,(i-1)+\beta
(j+1))/\alpha}
\end{equation}
(this computation is omitted).
Finally, taking into consideration the additivity of the Poincar\'e-Pontryagin-Melnikov function (\ref{m1})
with respect to the monomials of $P(x,y,0), Q(x,y,0)$
we conclude
$$ I(x) \, = \, x^{1-(\omega+\beta)/\alpha} \, \sum_{(i,j) \in \mathcal{J}}
c_{ij}\, x^{(\alpha\,i+\beta j)/\alpha},
$$
where $\mathcal{J} \, = \, \left\{ (i,j) \in \mathbb{Z}^2 \, : \,
-1 \leq i \leq d,\,  0 \leq j  \leq d+1, \, 0 \leq i+j \leq d
\right\}$ and $c_{ij} \in \mathbb{R}$.
 The
number $(d+1)(d+4)/2$ is the cardinal of the set $\mathcal{J}$. Obviously the function $I(x)$ has at most $(d+1)(d+4)/2-1$ zeros (counted with multiplicity), on the interval $(0,\infty)$.
This completes the proof of Theorem \ref{theo}.\bbox

\begin{remark}Clearly, the above bound is exact if and only if $\alpha$ and $\beta$ are not commensurable. The
weights $\alpha$ and $\beta$ of any weight--homogeneous polynomial
system {\rm (\ref{eqg})} with a center at the origin are
necessarily commensurable. To show this fact just take any
monomial of $f(x,y)$ with a nonzero coefficient, $a_{ij} x^i y^j$,
and any monomial of $g(x,y)$ with a nonzero coefficient: $b_{i'j'}
x^{i'} y^{j'}$. Since $f(\rho^\alpha\, x, \rho^\beta\, y) \, = \,
\rho^{\omega+\alpha} \, f(x,y)$ and $g(\rho^\alpha\, x,
\rho^\beta\, y) \, = \, \rho^{\omega+\beta} \, g(x,y)$, for any
$x,y,\rho \in \mathbb{R}$, we deduce the identities: $\alpha i +
\beta j \, = \, \omega + \alpha$ and $\alpha i' + \beta j' \, = \,
\omega + \beta$. We subtract them to deduce that $\alpha (i-i'-1)
\, + \, \beta (j-j'+1) \, = \, 0$, which gives that $\alpha$ and
$\beta$ are commensurable. In the case when both $f(x,y)$ and
$g(x,y)$  have only one monomial with nonzero coefficient
($f(x,y)=a_{ij} x^{i}y^j$, $g(x,y)=b_{i'j'} x^{i'} y^{j'}$) and
such that $i=i'+1$ and $j'=j+1$ the weights $\alpha$ and $\beta$
are not commensurable, but the origin is a linear node instead of
a center.
\end{remark}
\begin{remark}
The same result directly follows from the formula {\rm (\ref{m1})}
for $I(h)$; it suffices to note that $H, H_x, H_y, V$ are
weight-homogeneous functions of appropriate degree. \label{rem5}
\end{remark}

{\em Proof of Theorem} \ref{prop3}. The weight--homogeneous degree
of the system $(\alpha, \beta, 2\delta-\alpha-\beta)$ follows from
straightforward computations. \par The condition that $h_2(x,y)
\geq 0$ implies that the variation of the argument of $h_1(x,y)+i
h_2(x,y)$ along any closed path $l\subset \mathbb{R}^2\setminus
\{(0,0)\}$ is zero. Therefore for every fixed $\sigma \in
\mathbb{C}$ the function $(h_1+i h_2)^\sigma$ has a single valued
analytic continuation on $\mathbb{R}^2\setminus \{(0,0)\}$. From
now on we fix some determination of $(h_1+i h_2)^\sigma$. We note
that functions $H,V$ defined above, as well the associated
differential system, do not depend on this particular
determination. We may suppose without loss of generality that $\Re
\sigma > 0$ (otherwise we just replace $H$ by $1/H$). Then $H$ has
a continuous limit at $(0,0)$ and we may put $H(0,0)=0$. We claim
that each level set $\{(x,y): H(x,y)=\varepsilon \}$, $\varepsilon
>0$, is a smooth closed curve containing the origin. Indeed, the
restriction of $H$ on a half-line $l$ starting at the origin is
again a positive weight--homogeneous function. It follows that
$\{(x,y): H(x,y)=\varepsilon \} \cap l$ consists of a single point
and therefore $\{(x,y): H(x,y)=\varepsilon \}$ is a closed curve.
Suppose that $dH(x_0,y_0)=0$. Then the differential of $H$ is zero
at any point belonging to the half line $l_0$ starting at the
origin and containing $(x_0,y_0)$. It follows that $H|l_0$ is a
constant and moreover this constant equals to $0=H(0,0)$ which is
impossible. Therefore the level set $\{(x,y): H(x,y)=\varepsilon
\}$, $\forall \varepsilon >0$, is a closed periodic orbit. The
system has a global center.
 \bbox \newline

{\em Proof of Theorem} \ref{th1}. According to Theorem \ref{prop3}  the origin of system (\ref{eq1}) is a
center. As noted in Remark \ref{rem5}, the first integral $H$ and the inverse integrating factor $V$ are
weight--homogeneous. More precisely
\begin{lemma}
Given any $\rho \in \mathbb{C}-\{0\}$ we have:
\[ V\left( \rho\, x, \, \rho^{2n} y \right) \, = \,
\rho^{4n} \, V(x,y), \qquad H\left( \rho \, x, \, \rho^{2n} y
\right) \, = \, \rho^{4n} \, H(x,y), \] for any $(x,y) \in
\mathbb{R}^2-\left\{(0,0) \right\}$. \label{lem1}
\end{lemma}
The proof is a straightforward computation.

We define the following functions:
\[ I_{ij}(h) \, = \, \oint_{H=h} \frac{x^i\, y^j}{V(x,y)} \, dx,
\qquad J_{ij}(h) \, = \, \oint_{H=h} \frac{x^i\, y^j}{V(x,y)} \, dy, \] where $i,j$ are non--negative integer
numbers. By the expression of $I(h)$ and taking into account that $P(x,y)$ and $Q(x,y)$ are polynomials with
real coefficients of degree at most $4n-1$, we have that: \begin{equation} I(h) \, = \, \sum_{0 \leq i+j \leq
4n-1} \alpha_{ij}\, I_{ij}(h) \, + \, \beta_{ij}\, J_{ij}(h), \label{eq3} \end{equation} with $\alpha_{ij}$ and
$\beta_{ij}$ real numbers and $i,j$ are nonnegative integer numbers.
\begin{lemma} Given any $\rho, h \in \mathbb{R}$ with $h>0$ and
$\rho \neq 0$, we have:
\[ I_{ij} \left( \rho^{4n} h \right) \, = \, |\rho|^{i+1+2nj-4n}
I_{ij}(h), \qquad J_{ij} \left(\rho^{4n} h \right) \, = \,
|\rho|^{i+2n(j+1)-4n} J_{ij}(h), \] where $|\rho|$ stands for the
absolute value of $\rho$. \label{lem2}
\end{lemma}
{\em Proof.} The change of variables $x \mapsto |\rho| x$, $y \mapsto \rho^{2n} y$ in the integrals $I_{ij}$ and
$J_{ij}$ (which preserves  the orientation of the oval $\{H=h\}$) implies
\begin{eqnarray*}
I_{ij}\left( \rho^{4n} h \right) & = & \oint_{H=\rho^{4n} h}
\frac{x^i \, y^j}{V(x,y)} \, dx \, = \, \oint_{H=h}
|\rho|^{i+1+2nj-4n} \, \frac{x^i \, y^j}{V(x,y)} \, dx
\vspace{0.3cm}
\\ & = & |\rho|^{i+1+2nj-4n} I_{ij}(h). \vspace{0.4cm} \\ J_{ij}
\left(\rho^{4n} h \right) & = & \oint_{H=\rho^{4n} h} \frac{x^i \,
y^j}{V(x,y)} \, dy \, = \, \oint_{H=h} |\rho|^{i+2n(j+1)-4n} \,
\frac{x^i \, y^j}{V(x,y)} \, dy \vspace{0.3cm} \\ & = &
|\rho|^{i+2n(j+1)-4n} J_{ij}(h).
\end{eqnarray*} \bbox \newline

$ $From the previous lemma we give the form of the functions
$I_{ij}(h)$ and $J_{ij}(h)$. The same argument given in the proof
of Theorem \ref{theo} holds: we choose $h=1$ and
$\rho=h^{\frac{1}{4n}}$ in the expressions given in Lemma
\ref{lem2}.
\begin{lemma}
Given any $h>0$ we have:
\[ I_{ij}(h) \, = \, h^{\frac{i+1+2nj-4n}{4n}} \, I_{ij}(1), \qquad J_{ij}(h) \, = \, h^{\frac{i+2n(j+1)-4n}{4n}} \,
J_{ij}(1). \] \label{lem3}
\end{lemma}

We get that $I_{ij}(h)$ and $J_{ij}(h)$ are monomials of $h$ up to
a fractional power. Hence, we only need to determine how many of
these functions are linearly independent, so as to know how many
zeroes can have $I(h)$. We first determine which of the functions
$I_{ij}(h)$ and $J_{ij}(h)$ are identically zero.
\begin{lemma}
If $\, i$ is odd, then $I_{ij}(h) \equiv 0$. If $\, i$ is even,
then $J_{ij}(h) \equiv 0$. \label{lem4} \end{lemma} {\em Proof.}
Let us consider the change of variables $x \mapsto -x$ and $y
\mapsto y$ in the integrals $I_{ij}$ and $J_{ij}$. We note that
this change of coordinates reverses the orientation of the oval
$\{H=h\}$. We denote the oval with reversed orientation by
$-\{H=h\}$. Moreover, from Lemma \ref{lem1} we have that this
change of coordinates leaves $V(x,y)$ and $H(x,y)$ invariant.
\begin{eqnarray*}
I_{ij}(h) & = & \oint_{H=h} \frac{x^i\, y^j}{V(x,y)}\, dx \, = \,
\oint_{-\{H=h\}} (-1)^{i+1} \, \frac{x^i\, y^j}{V(x,y)}\, dx \, =
\, (-1)^i \, I_{ij}(h). \vspace{0.3cm} \\ J_{ij}(h) & = &
\oint_{H=h} \frac{x^i\, y^j}{V(x,y)}\, dy \, = \, \oint_{-\{H=h\}}
(-1)^{i} \, \frac{x^i\, y^j}{V(x,y)}\, dy \, = \, (-1)^{i+1} \,
J_{ij}(h).
\end{eqnarray*}
Taking $i$ odd, we deduce that $I_{ij}(h)$ needs to be identically
zero, and the same is true for $J_{ij}(h)$ taking $i$ even. \bbox
\par
Lemma \ref{lem4} can also be proved using Green's Theorem and
analogous reasonings. \newline

We are going to characterize some of the functions $I_{ij}(h)$ and
$J_{ij}(h)$ which are not zero at any point of $h>0$.
\begin{lemma}
For any $k, \ell$ non--negative integers we have:
\[ I_{2k,2\ell+1}(h) \, \not\equiv \, 0, \qquad J_{2k+1,2\ell}(h)
\, \not\equiv \, 0, \] where we recall that
\[ I_{2k,2\ell+1}(h) \, : = \, \oint_{H=h} \frac{x^{2k}
\, y^{2 \ell +1}}{V(x,y)} \, dx, \qquad J_{2k+1,2\ell}(h) \, := \,
\oint_{H=h} \frac{x^{2k+1} \, y^{2 \ell}}{V(x,y)} \, dy,\] with
$H(x,y)$ as given in {\rm (\ref{Hn})} and $V(x,y)\, := \,
y^2+x^{2n}y-c\, x^{4n}$. \label{lem5}
\end{lemma} {\em Proof.} We note that the orientation of system (\ref{eq1})
over the oval $H=h$ is clockwise and that $V(x,y)$ is strictly
positive over all the oval $H=h$. \par Let us denote by $x_1(h)$
and $x_2(h)$ the intersections of the oval $H=h$ with the
horizontal axis $(y=0)$ with $x_1(h)<0$ and $x_2(h)>0$. We denote
by $\{H=h\}_{y<0}$ the half part of the oval below the horizontal
axis oriented from $x_2(h)$ to $x_1(h)$ and by $\{H=h\}_{y>0}$ the
half part of the oval above the horizontal axis oriented from
$x_1(h)$ to $x_2(h)$. We have:
\begin{eqnarray*} I_{2k,2\ell+1}(h) & = & \oint_{H=h} \frac{x^{2k}
\, y^{2 \ell +1}}{V(x,y)} \, dx \, = \, \int_{x_2(h)}^{x_1(h)}
\frac{x^{2k} \, y^{2 \ell +1}}{V(x,y)} \, dx \, + \,
\int_{x_1(h)}^{x_2(h)} \frac{x^{2k} \, y^{2 \ell +1}}{V(x,y)} \,
dx \vspace{0.3cm} \\ & = & -\int_{x_1(h)}^{x_2(h)}  \frac{x^{2k}
\, y^{2 \ell +1}}{V(x,y)} \, dx \, + \, \int_{x_1(h)}^{x_2(h)}
\frac{x^{2k} \, y^{2 \ell +1}}{V(x,y)} \, dx,
\end{eqnarray*}
where the first integral is done over the path $\{H=h\}_{y<0}$ and
the second integral is done over the path $\{H=h\}_{y>0}$.
Therefore, the first integral is strictly negative and the second
integral is strictly positive and we are adding two positive
values, due to the minus sign. Hence, $I_{2k,2\ell+1}(h)>0$ and it
cannot be zero at any point. \par Analogously for
$J_{2k+1,2\ell}(h)$, we define $y_1(h)$ and $y_2(h)$ the
intersections of the oval $H=h$ with the vertical axis $(x=0)$
with $y_1(h)<0$ and $y_2(h)>0$. We denote by $\{H=h\}_{x<0}$ the
half part of the oval at the left of the vertical axis oriented
from $y_1(h)$ to $y_2(h)$ and by $\{H=h\}_{x>0}$ the half part of
the oval at the right of the vertical axis oriented from $y_2(h)$
to $y_1(h)$, and we have:
\begin{eqnarray*} J_{2k+1,2\ell}(h) & = & \oint_{H=h} \frac{x^{2k+1}
\, y^{2 \ell}}{V(x,y)} \, dy \, = \, \int_{y_2(h)}^{y_1(h)}
\frac{x^{2k+1} \, y^{2 \ell}}{V(x,y)} \, dy \, + \,
\int_{y_1(h)}^{y_2(h)} \frac{x^{2k+1} \, y^{2 \ell}}{V(x,y)} \, dy
\vspace{0.3cm} \\ & = & -\int_{y_1(h)}^{y_2(h)}  \frac{x^{2k+1} \,
y^{2 \ell}}{V(x,y)} \, dy \, + \, \int_{y_1(h)}^{y_2(h)}
\frac{x^{2k+1} \, y^{2 \ell}}{V(x,y)} \, dy,
\end{eqnarray*}
where the first integral is done over the path $\{H=h\}_{x>0}$ and
the second integral is done over the path $\{H=h\}_{x<0}$.
Therefore, the first integral is strictly positive and the second
integral is strictly negative and we are adding two negative
values, due to the minus sign. Hence $J_{2k+1,2\ell}(h)<0$ and it
cannot be zero at any point. \bbox \newline

In fact, we have been able to numerically prove that the integrals
$I_{2k,2\ell}(h) \, \not\equiv \, 0$ and $J_{2k+1,2\ell+1}(h) \,
\not\equiv \, 0$ for some particular fixed values of the integers
$k,\ell$, with $k\geq 0$ and $\ell\geq 0$, in some fixed cases of
the function $H(x,y)$ defined in (\ref{Hn}). To show this fact, we
have parameterized the oval $H(x,y)=h$ by $(x_{+}(\tau),y(\tau))$
when $x>0$ and $(x_{-}(\tau),y(\tau))$ when $x<0$, with
\[ x_{\pm}(\tau) \, = \, \pm \, h^{\frac{1}{4n}} \,
(\tau+\mu_{+})^{-\frac{\sigma}{4n}} \,
(\tau+\mu_{-})^{-\frac{\bar{\sigma}}{4n}}, \quad  y(\tau) \, = \,
h^{\frac{1}{2}} \, \tau \, (\tau+\mu_{+})^{-\frac{\sigma}{2}} \,
(\tau+\mu_{-})^{-\frac{\bar{\sigma}}{2}}, \] where $\mu_{\pm} \, =
\, (1\pm\sqrt{1+4c})/2$, $\sigma=1-1/\sqrt{1+4c}$ and the rank of
the parameter $\tau$ is all the real line $\tau \in
(-\infty,+\infty)$ in both parts of the oval. When we write
$h^{\frac{1}{4n}}$ or $h^{\frac{1}{2}}$ we mean the positive real
root. \newline

Next lemma shows that the possible nonzero values of the integrals
$J_{ij}(h)$ are redundant, since there is an integral of the form
$I_{ij}(h)$ which corresponds to the same monomial.
\begin{lemma} The first Melnikov function can be expressed by:
\[ I(h) \, = \, \sum_{(i,j) \in \mathcal{I}} \alpha_{ij} \,
I_{ij}(h), \] where $\alpha_{ij} \in \mathbb{R}$ and $\mathcal{I}
\, = \, \left\{ (i,j)  \in  \mathbb{Z}^2 \, : \, 0\leq i,j \leq
4n-1, \, i \, \, \mbox{is even}, \, i+j\leq 4n-1\right\}$.
\label{lem6}
\end{lemma} {\em Proof.} The expression of $I(h)$ given in
(\ref{eq3}) ensures that this function is a linear combination of the integrals $I_{ij}(h)$ and $J_{ij}(h)$. We
are going to show that any possible monomial expressed by a $J_{ij}(h)$ can also be got by a monomial of
$I_{i'j'}(h)$. Lemma \ref{lem4} gives that only the integrals $J_{ij}(h)$ with $i$ odd need to be considered,
hence we take any two nonnegative integers $(i,j)$ such that $i$ is odd and $i+j\leq 4n-1$. We define $i'=i-1$
and $j'=j+1$ and we have that both $i'$ and $j'$ are nonnegative integers strictly lower than $4n$ and
$i'+j'=i+j\leq 4n-1$. Moreover, $i'$ is even in accordance with Lemma \ref{lem4} and $(i',j') \in \mathcal{I}$.
Hence, any monomial given by a $J_{ij}(h)$ is also expressed by a $I_{ij}(h)$ with $i$ even. \bbox \newline

To end with the proof of Theorem \ref{th1}, we first count the
cardinal of $\mathcal{I}$ and we are going to show that $\sharp \,
\mathcal{I}\, = \, 2n(2n+1)$. We take any even nonnegative integer
$i$ from $0$ to $4n-2$, that is, we take $2n$ possible values of
$k$ with $i=2k$ and given a fixed $0\leq k\leq 2n-1$ we can take
any value of $j$ from $0$ to $4n-1-2k$. Hence,
\[ \sharp \, \mathcal{I} \, = \, \sum_{k=0}^{2n-1} \,
\sum_{j=0}^{4n-1-2k} \, 1 \, = \, \sum_{k=0}^{2n-1} \, (4n-2k) \,
= \, 2 \, + \, 4 \, + \, 6 \, + \, \ldots \, + \, 4n \, = \,
2n(2n+1).\] We note that the cardinal of $\mathcal{I}$ is an upper
bound for the number of independent monomials given by the nonzero
$I_{ij}(h)$ because it may happen that two elements of $I_{ij}(h)$
give rise to the same monomial, that is, it may happen that $(i,j)
\in \mathcal{I}$, $(i',j') \in \mathcal{I}$ with $(i,j) \neq
(i',j')$ but $i+2nj=i'+2nj'$. To end with the proof of part {\rm
(a)} of Theorem \ref{th1}, we need to characterize the number of
repeated exponents corresponding to different indexes $(i,j) \in
\mathcal{I}$, $(i',j') \in \mathcal{I}$. We note that if two such
pairs give $i+2nj=i'+2nj'$, then $j$ and $j'$ can differ at most
by one, because if they differ by two or more then $i$ and $i'$
differ by $4n$ or more which is not possible since $0 \leq
i,i'\leq 4n-1$. We assume that $j'=j+1$ and we have that
$i'=i-2n$. Hence, fixing $(i,j)$, the condition to have a repeated
exponent is that $i \geq 2n$ and we already have that $j<4n-1$
because $i+j\leq 4n-1$. Let us count how many of these indexes we
do have in $\mathcal{I}$. We fix $k$ such that $i=2k$ and $k$ goes
from $n$ to $2n-1$. Given such a $k$ we can take, as before, any
value of $j$ from $0$ to $4n-1-2k$. Therefore, the number of
repeated exponents is:
\[ \sum_{k=n}^{2n-1} \sum_{j=0}^{4n-1-2k} \, 1 \, = \,
\sum_{k=n}^{2n-1}(4n-2k) \, = \, 2 \, + \, 4 \,  + \, \ldots \, +
\, 2n \, = \, n(n+1). \] We conclude that the number of different
exponents associated to the indexes of $\mathcal{I}$ is
$2n(2n+1)-n(n+1) \, = \, n(3n+1)$. \par We have given an upper
bound for the number of independent functions in which $I(h)$ can
be split as a linear combination. Since this upper bound is
$n(3n+1)$, we have that an upper bound for the number of isolated
zeroes of $I(h)$ is $n(3n+1)-1$.
\par The proof of part (b) in Theorem \ref{th1} comes from Lemma
\ref{lem5}, which ensures that only half of the $2n(2n+1)$ functions $I_{ij}(h)$ are ensured to be different
from zero, the ones with an even $i$ and an odd $j$. To end with, we only need to see that  two such functions
give place to two independent monomials. Let $(2k,2\ell+1) \in \mathcal{I}$ and $(2k',2\ell'+1) \in \mathcal{I}$
with $(2k,2\ell+1)  \neq (2k',2\ell'+1)$ be such that $2k+2n(2\ell+1)=2k'+2n(2\ell'+1)$ and we will get a
contradiction. If $k=k'$ then $\ell=\ell'$, so we can assume that $k<k'$. We have $2n(\ell-\ell')=k'-k$ which
gives that $\ell \geq \ell'+1$. This inequality gives that $k' \geq  2n+k$ and since $k \geq 0$, we have that
$k' \geq 2n$, which is impossible since $2k'\leq 4n-1$. Hence, any two exponents given by different functions
$I_{2k, 2 \ell+1}(h)$ are independent. Choosing adequate parameters in system (\ref{eq1}) we get that $I(h)$ can
always be given as a sum of at least $n(2n+1)$ independent functions and, adapting parameters, it can always be
chosen with at least $n(2n+1)-1$ different isolated zeroes. \bbox

\appendix
\section{The Andreev theorem}

 We already noted that Theorem \ref{prop3} implies that
the origin of the system (\ref{eq1}) is a global center. It is
interesting to note that the same result follows also from a
classical theorem of Andreev\cite{Andreev} which we briefly
illustrate. This approach can be useful in other cases of
nilpotent centers, which are not covered by Theorem \ref{prop3}.
\begin{theorem}[Andreev's Theorem \cite{Andreev}] Let $F(x,y)$ and
$G(x,y)$ be analytic functions in a neighborhood of the origin of order $\geq 2$ and let the origin be an
isolated singularity of the differential system: $\dot{x}=y+F(x,y)$, $\dot{y}=G(x,y)$. Let $y=\phi(x)$ be the
solution of the equation $y+F(x,y)=0$ such that $\phi(0)=0$. We denote by $\xi(x)=G(x,\phi(x))$ and $\Delta(x) =
{\rm div}(x, \phi(x))$ and we develop them in a neighborhood of $x=0$:
\[ \xi(x) \, = \, \alpha_1 \, x^{k_1} \, + \,
\mathcal{O}\left(x^{k_1+1}\right), \qquad \Delta(x) \, = \, \alpha_2 \, x^{k_2} \, + \,
\mathcal{O}\left(x^{k_2+1}\right), \] where $\alpha_1 \neq 0$, $k_1 \geq 2$ and $\alpha_2 \neq 0$, $k_2\geq 1$
or $\Delta(x) \equiv 0$. The origin of the differential system is monodromic if, and only if, $\alpha_1<0$,
$k_1$ is an odd number, and one of the following three conditions holds: \begin{enumerate}
\item[{\rm (a)}] $k_1 \, = \, 2 k_2 + 1$ and $\alpha_2^2+4\alpha_1(k_2+1)<0$,
\item[{\rm (b)}] $k_1<2k_2+1$,
\item[{\rm (c)}] $\Delta(x) \equiv 0$. \end{enumerate} \label{thandreev}
\end{theorem}

The function $\xi(x)$ associated to  (\ref{eq1}) is $\xi(x)=2ncx^{4n-1}$, so $\alpha_1=2nc$ and $k_1=4n-1$. The
divergence of system (\ref{eq1}) is ${\rm div}(x,y)=2nx^{2n-1}$, so $\Delta(x)=2nx^{2n-1}$ and we have
$\alpha_2=2n$ and $k_2=2n-1$. Since $n$ is a natural number and $c<-1/4$, we already have that $\alpha_1<0$ and
$k_1$ is an odd integer. The condition {\rm (a)} of Theorem \ref{thandreev} is satisfied because $2 k_2 +1 =
4n-1=k_1$ and $\alpha_2^2+4\alpha_1(k_2+1)=4n^2(1+4c)<0$. We conclude that the origin of system (\ref{eq1}) is
monodromic. Finally, the origin of the system (\ref{eq1}) is time--reversible because it is invariant by the
change $(x,y,t) \to (-x,y,-t)$, which implies that it is a center. Moreover the center is global, as the system
is weight--homogeneous.

\vspace{2ex} \noindent {\bf Acknowledgment.}  We thank I.D. Iliev for the critical remarks.


\begin{thebibliography}{99}

\bibitem{Andreev} {\sc A. Andreev},{\it \ Investigations of the behavior
of the integral curves of a system of two differential equations
in the neighborhood of a singular point.} Translation Amer. Math.
Soc., {\bf 8} (1958), 187--207.

\bibitem{Francoise} {\sc J.-P. Fran\c{c}oise},{\it \ Successive
derivatives of a first return map, application to the study of
quadratic vector fields.} Ergod. Theory Dynam. Syst. {\bf 16}
(1996), 87--96.

\bibitem{gav} {\sc L. Gavrilov},
{\it \ Higher order Poincare-Pontryagin functions and iterated
path integrals.}  Ann. Fac. Sci. Toulouse Math. (6) {\bf 14}
(2005), no. 4, 663--682.

%


\bibitem{Iliev} {\sc I.D. Iliev},{\it \ On second order
bifurcations of limit cycles.} J. Lond. Math. Soc. {\bf 58}
(1998), 353--366.

\bibitem{Ilya} {\sc Yu. S. Ilyashenko},{\it \  The appearance of limit cycles under a perturbation of the equation
$dw/dz=-R\sb{z}/R\sb{w}$, where $R(z,\,w)$ is a polynomial.}
(Russian) Mat. Sb. (N.S.) {\bf 78 (120)} (1969), 360--373.

\bibitem{Li} {\sc C. Li, W. Li, J. Llibre and Z. Zhang}, {\it \ Polynomial systems: a lower bound for the weakened
16th Hilbert problem.} Extracta Math. {\bf 16} (2001), no. 3,
441--447.

\bibitem{LLYZ} {\sc W. Li, J. Llibre, J. Yang and Z. Zhang}, {\it \ Limit cycles
bifurcating from the period annulus of quasi-homogeneous centers.}
To appear in J. Dyn. Diff. Equat.

\bibitem{Llibre} {\sc J. Llibre and X. Zhang},{\it \ On the number of limit cycles for
some perturbed Hamiltonian polynomial systems.}  Dyn. Contin.
Discrete Impuls. Syst. Ser. A Math. Anal. {\bf 8} (2001), no. 2,
161--181.

\bibitem{Perko} {\sc L. Perko}, {\it \ Differential equations and dynamical
systems.} Third edition. Texts in Applied Mathematics, {\bf 7}.
Springer--Verlag, New York, 2001.

\bibitem{Spivak} {\sc M. Spivak}, {\it \ A comprehensive introduction to
differential geometry. Vol. One.} Published by M. Spivak, Brandeis
Univ., Waltham, Mass. 1970.

\bibitem{ZhaoZhang} {\sc Y. Zhao and Z. Zhang},{\it \
Abelian integrals and period functions for quasihomogeneous
Hamiltonian vector fields.} Appl. Math. Lett. {\bf 18} (2005),
563--569.

\end{thebibliography}
\end{document}